\documentclass{gtart_h}

\def\ifplaintex{\expandafter\ifx\csname documentclass\endcsname\relax}

\def\gtp{{\mathsurround=0pt\it $\cal G\mskip-2mu$eometry \&\ 
$\cal T\!\!$opology $\cal P\!$ublications}}  

\def\recd{{\small Received:\qua\receiveddate\ifx\reviseddate\relax
\else\qquad Revised:\qua\reviseddate\fi\par}} 

\def\lognumber#1{\def\thelognumber{#1}}
\def\volumenumber#1{\def\thevolumenumber{#1}}
\def\volumeyear#1{\def\thevolumeyear{#1}}
\def\papernumber#1{\def\thepapernumber{#1}}
\def\pagenumbers#1#2{\def\startpage{#1}\def\finishpage{#2}}
\def\published#1{\def\publishdate{#1}}

\def\received#1{\def\receiveddate{#1}}

\def\accepted#1{\def\accepteddate{#1}}

\def\asciiemail#1{\def\theasciiemail{#1}}

\let\\\par\let\thelognumber\relax\let\thevolumenumber\relax
\let\thepapernumber\relax\let\thevolumeyear\relax\let\startpage\relax
\let\finishpage\relax\let\publishdate\relax\let\receiveddate\relax
\let\reviseddate\relax\let\accepteddate\relax\let\theasciititle\relax
\let\theasciiauthors\relax
\let\theasciiabstract\relax

\let\theasciiemail\relax

\ifplaintex
\font\logobig=cmssbx10 scaled 3836
\font\logomed=cmssbx10 scaled 2557
\else
\font\logobig=cmssbx10 scaled 4200
\font\logomed=cmssbx10 scaled 2800
\fi

\long\def\makeagttitle{   
\count0=\startpage
\agt\hfill      
\hbox to 45truept{\vbox to 0pt{\vglue -13truept{\logomed A\kern -.37em{\logobig 
T}\kern -.38em G}\vss}\hss}
\break
{\small Volume \thevolumenumber\ (\thevolumeyear)
\startpage--\finishpage\nl
Published: \publishdate}

\vglue .25truein

{\parskip=0pt\leftskip 0pt plus
1fil\def\\{\par\smallskip}{\Large\bf\thetitle}\par\medskip} \vglue
0.05truein

{\parskip=0pt\leftskip 0pt plus 1fil\def\\{\par}{\sc\theauthors}
\par\medskip}%
 
\vglue 0.03truein 

{\small\leftskip 25truept\rightskip 25truept{\bf Abstract}\stdspace\theabstract

{\bf AMS Classification}\stdspace\theprimaryclass
\ifx\thesecondaryclass\relax\else; \thesecondaryclass\fi\par
{\bf Keywords}\stdspace \thekeywords\par}\vglue 7truept

}   

\ifplaintex
\font\phead=cmsl9 scaled 950
\font\pnum=cmbx10 scaled 913
\font\pfoot=cmsl9 scaled 950
\headline{\vbox to 0pt{\vskip -4.5mm\line{\small\phead\ifnum
\count0=\startpage ISSN 1472-2739 (on-line) 1472-2747 (printed)
\hfill {\pnum\folio}\else\ifodd\count0\def\\{ }%
\ifx\theshorttitle\relax\thetitle\else\theshorttitle\fi\hfill{\pnum\folio}
\else\def\\{ and }{\pnum\folio}\hfill\ifx\theshortauthors\relax\theauthors
\else\theshortauthors\fi\fi\fi}\vss}}
\footline{\vbox to 0pt{\vglue 0mm\line{\small\pfoot\ifnum\count0=\startpage
\copyright\ \gtp\hfill\else
\agt, Volume \thevolumenumber\ (\thevolumeyear)\hfill\fi}\vss}}
\else
\font=cmsl9 scaled 1050
\font\lnum=cmbx10 
\font=cmsl9 scaled 1050
\makeatletter
\def\@oddhead{{\small\lhead\ifnum\count0=\startpage ISSN 1472-2739 
(on-line) 1472-2747 (printed)\hfill {\lnum\number\count0}\else\ifodd\count0
\def\\{ }\ifx\theshorttitle\relax \thetitle \else\theshorttitle\fi\hfill
{\lnum\number\count0}\else\def\\{ and }{\lnum\number\count0}
\hfill\ifx\theshortauthors\relax 
\theauthors\else\theshortauthors\fi\fi\fi}}\def\@evenhead{\@oddhead}
\def\@oddfoot{\small\lfoot\ifnum\count0=\startpage\copyright\ \gtp\hfill\else
\agt, Volume \thevolumenumber\ (\thevolumeyear)\hfill\fi}
\def\@evenfoot{\@oddfoot}
\makeatother
\fi
\let\maketitlepage\makeagttitle
\let\makeshorttitle\maketitlepage
\let\maketitle\maketitlepage

\newwrite\gtoutfile
\long\gdef\makeheadfile{  
{\def\\{, }\def\s{ }
\immediate\openout\gtoutfile head.xxx
\immediate\write\gtoutfile{Proxy-for: \ifx\theasciiauthors\relax
\theauthors\else\theasciiauthors\fi\s<\ifx\theasciiemail\relax\theemail\else\theasciiemail\fi>}
\immediate\write\gtoutfile{\noexpand\\}
\immediate\write\gtoutfile{Authors: \ifx\theasciiauthors\relax
\theauthors\else\theasciiauthors\fi}
{\def\\{ }\immediate\write\gtoutfile{Title: \ifx\theasciititle\relax
\thetitle\else\theasciititle\fi}}
\immediate\write\gtoutfile{Subj-class: GT or SG, GR etc}
\immediate\write\gtoutfile{MSC-class: \theprimaryclass\ifx\thesecondaryclass\relax\else, \thesecondaryclass\fi}
\immediate\write\gtoutfile{Journal-ref: Algebr. Geom. Topol. \thevolumenumber\s
(\thevolumeyear) \startpage-\finishpage}
\immediate\write\gtoutfile{Comments: Published by Algebraic and
Geometric Topology at}
\immediate\write\gtoutfile{\s\s\s  http://www.maths.warwick.ac.uk/agt/AGTVol\thevolumenumber/agt-\thevolumenumber-\thepapernumber.abs.html}
\immediate\write\gtoutfile{\noexpand\\}
\immediate\write\gtoutfile{}
\ifx\theasciiabstract\relax
\immediate\write\gtoutfile{\theabstract}\else
\immediate\write\gtoutfile{\theasciiabstract}\fi
\immediate\write\gtoutfile{}
\immediate\write\gtoutfile{\noexpand\\}
\immediate\write\gtoutfile{}
\immediate\closeout\gtoutfile}}  

\def\maketitlepage{\makeagttitle\makeheadfile}
\let\makeshorttitle\maketitlepage
\let\maketitle\maketitlepage

\lognumber{53}
\volumenumber{5}
\volumeyear{2005}
\papernumber{53}
\pagenumbers{1365}{1388}
\received{16 December 2004} 
\accepted{11 May 2005}
\published{14 October 2005}

\usepackage{graphicx,amsmath,amssymb}
\usepackage[mathscr]{eucal}

\def\figref#1{\hyperlink{#1anchor}{Figure~\ref*{#1}}}
\def\anchor#1{\noindent\hypertarget{#1anchor}{\smash{$\phantom{99}$}}}

\newcommand{\pgl}{\mbox{$P_G(\lambda)$}}

\newcommand\M{{\mathcal M}}
\newcommand\N{{\mathcal N}}
\newcommand\Z{{\mathbb{Z}}}

\newtheorem{theorem}{Theorem}[section]
\newtheorem{corollary}[theorem]{Corollary}
\newtheorem{lemma}[theorem]{Lemma}
\newtheorem{proposition}[theorem]{Proposition}

\theoremstyle{remark}
\newtheorem{definition}[theorem]{Definition}
\newtheorem{example}[theorem]{Example}

\newtheorem{remark}[theorem]{Remark}

\begin{document}

\title{A categorification for the chromatic polynomial}
\author{Laure Helme-Guizon\\Yongwu Rong}

\address{Department of Mathematics, The George Washington
University\\Washington, DC 20052, USA}
\gtemail{\mailto{lhelmeg@yahoo.com}, \mailto{rong@gwu.edu}}
\asciiemail{lhelmeg@yahoo.com, rong@gwu.edu}

\begin{abstract}
For each graph we construct graded cohomology groups whose graded
Euler characteristic is the chromatic polynomial of the graph. We
show the cohomology groups satisfy a long exact sequence which
corresponds to the well-known deletion-contraction rule. This work
is motivated by Khovanov's work on categorification of the Jones
polynomial of knots.
\end{abstract}



\maketitle
\section{Introduction}

In recent years, there have been a great deal of interests in
Khovanov cohomology theory introduced in \cite{K00}. For each link
$L$ in $S^{3}$, Khovanov defines a family of cohomology groups
$H^{ij}(L)$ whose Euler characteristic $\sum_{i,j}(-1)^{i}  q^j rk
(H^{i,j}(L)) $ is the Jones polynomial of $L$. These groups were
constructed through a categorification process which starts with a
state sum of the Jones polynomial, constructs a group for each
term in the summation, and then defines boundary maps between
these groups appropriately. More recently, Khovanov and Rozansky
have extended the theory for the HOMFLYPT polynomial \cite{KR04}
\cite{KR05}.

It is natural to ask if similar categorifications can be done for
other invariants with state sums. In this paper, we establish a
cohomology theory that categorifies the chromatic polynomial for
graphs. This theory is based on the polynomial algebra with one
variable $x$ satisfying $x^2=0$. We show our cohomology theory
satisfies a long exact sequence which can be considered as a
categorification for the well-known deletion-contraction rule for
the chromatic polynomial. This exact sequence helps us to compute
the cohomology groups of several classes of graphs. In particular,
we point out that torsions do occur in the cohomology for some
graphs.

In Section 2, we show our construction. This is explained in two
equivalent settings: the cubic complex approach and the enhanced
state approach. In Section 3, we discuss various properties of
these cohomology groups, including the long exact sequence.
Several computational examples are given. Some questions and more
recent developments are discussed in Section 4.


\medskip \textbf{Acknowledgements}\qua We wish to thank Joe
Bonin for sharing his knowledge about graphs, and Mikhail Khovanov
for his comments. The first author was partially supported by the
CCAS Dean's Fellowship at the George Washington University. The
second author was partially supported by NSF grant DMS-0513918.

\section{From a graph to cohomology groups}

\subsection{A diagram for the chromatic polynomial}
We begin with a brief review for the chromatic polynomial.
 Let $G$ be a graph with vertex set $V(G)$ and edge set $E(G)$.
For each positive integer $\lambda$, let $\{ 1, 2, \cdots,
\lambda\}$ be the set of $\lambda$-colors. A
\emph{$\lambda$-coloring} of $G$ is an assignment of a
$\lambda$-color to vertices of $G$ such that vertices that are
connected by an edge in $G$ always have different colors. Let
$P_G(\lambda)$ be the number of $\lambda$-colorings of $G$. It is
well-known that  $P_G(\lambda)$ satisfies the
\emph{deletion-contraction relation}

\[ P_G(\lambda)=P_{G-e}(\lambda) - P_{G/e}(\lambda) \]
Furthermore, it is obvious that
\[ P_{N_v}(\lambda)=\lambda^v \mbox{ where $N_v$ is the graph with $v$
vertices and no edges.} \]

These two equations uniquely determines $P_G(\lambda)$. They also
imply that $P_G(\lambda)$ is always a polynomial of $\lambda$, known
as the {\it chromatic polynomial}.

There is another formula for \pgl\ that is useful for us. For each
$s\subseteq E(G)$, let $[G:s]$ be the graph whose vertex set is
$V(G)$ and whose edge set is $s$, let $k(s)$ be the number of
connected components of $[G:s]$. We have
\begin{equation} \label{State Sum 1}
 P_G(\lambda)=\mathrel{\mathop{\sum }\limits_{s\subseteq
E(G)}}(-1)^{\left| s\right| }\lambda ^{k(s)}.
\end{equation}
Equivalently, grouping the terms $s$ with the same number of edges
yields state sum formula
\begin{equation}
\label{State sum Chrom} P_G(\lambda)=\mathrel{\mathop{\sum
}\limits_{i\geq 0}}(-1)^{i}\mathrel{\mathop{\sum
}\limits_{s\subseteq E(G),\left| s\right| =i}}\lambda ^{k(s)}.
\end{equation}
Formula (\ref{State Sum 1}) follows easily from the well-known
state sum formula for the Tutte polynomial
$T(G,x,y)=\mathrel{\mathop{\sum }\limits_{s\subseteq
E(G)}}(x-1)^{r(E)-r(s)}(y-1)^{\left| s\right| -r(s)},$ where the
rank function $r(s)$ is the number of vertices of $G$ minus the
number of connected components of $[G:s]$. One simply applies the
known relation
$P(G,\lambda)=(-1)^{r(E)}\lambda^{k(G)}T(G;1-\lambda,0)$ between
the two polynomials.

Our chain complex will depend on an ordering of the edges of the
graph. Thus, let $G=(V,E)$ be a graph with an ordering on its $n$
edges. For each $s\subseteq E$, the spanning subgraph $[G:s]$ (a
\textit{spanning subgraph} of $G$ is one that contains all the
vertices of $G$) can be described unambiguously by an element
$\varepsilon =(\varepsilon _{1},\varepsilon _{2},....,\varepsilon
_{n})$ of $\left\{ 0,1\right\} ^{n}$ with the convention
$\varepsilon _{k}=1$ if the edge $e_{k}$ is in $s$ and
$\varepsilon _{k}=0$ otherwise. This $ \varepsilon $ is called the
\textit{label} of $s$ and will be denoted by $ \varepsilon _{s}$
or simply $\varepsilon .$ Conversely, to any $\varepsilon \in $
$\left\{ 0,1\right\} ^{n}$,$\ $we can associate a set
$s_{\varepsilon } $ of edges of $G$ that corresponds to
$\varepsilon .$ When we  think of $s$ in terms of the label
$\varepsilon$, we may refer to the graph
$\ [G:s]$ as $ G_{\varepsilon }.$ 

\begin{figure}[ht!]\anchor{figure State sum Chrom}
\begin{center}
\cl{\scalebox{.6}{\includegraphics{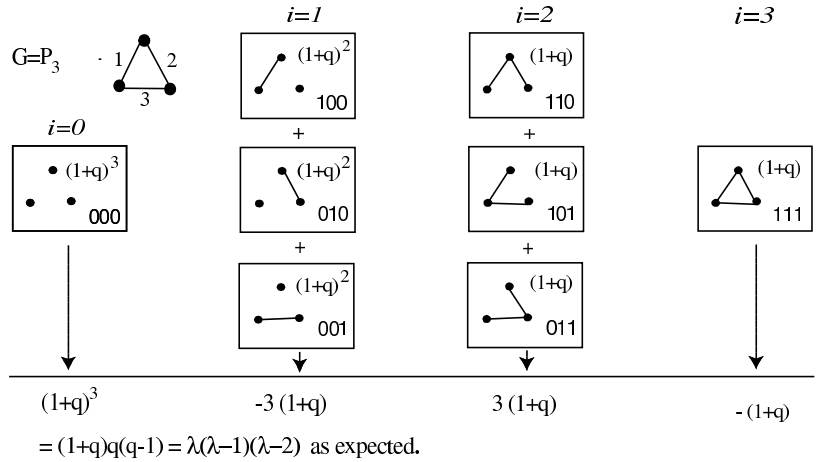}}}
\end{center}\vglue-.5cm
\caption{Diagram for a state sum computation of the
Chromatic Polynomial when $G=P_3$} \label{figure State sum Chrom}
\end{figure}

The state sum formula in (\ref{State sum Chrom}) can be
illustrated in \figref{figure State sum Chrom}. With the
variable change $\lambda =1+q$, this gives the chromatic
polynomial.

Each subset of edges $s$ (represented in \figref{figure State
sum Chrom} by a labelled rectangle) corresponds to a term in the
state sum and therefore will be called a \textit{state}.
Equivalently, if we think of a state in term of its label, we may
call it a vertex of the cube $\{0,1\}^{n}$.

Each state corresponds to a subset $s$ of $E,$ the n-list of $0^{\prime }s$
and $1^{\prime }s$ at the bottom of each rectangle is its label and the term
of the form $(1+q)^{k(s)}=\lambda ^{k(s)}$ is its contribution to the
chromatic polynomial (without sign).

Note that we have drawn all the states that have the same number of
edges in the same column, so that the column with label $i=i_{0}$
contains all the states with $i_{0}$ edges. Such states are called
the states of \textit{height} $i_{0}$. We denote the height of a
state with label $\varepsilon $ by $\left| \varepsilon \right|$.

\subsection{The cubic complex construction of the chain complex}

\subsubsection{The chain groups}

We are going to assign a graded $\mathbb{Z}$-module to each state
and define a notion of graded dimension so that the $(1+q)^{k}$
that appears in the rectangle is the graded dimension of this
$\mathbb{Z}$-module. The construction is similar to Bar-Natan's
description for the Khovanov cohomology for knots and links \cite
{BN02}.

\paragraph{Graded dimension of a graded $\mathbb{Z}$-module}

\begin{definition}
Let $\M=\mathrel{\mathop{\oplus }\limits_{j}}M_{j}$ be a graded
$\mathbb{Z}$-module where $\left\{ M_{j}\right\} $ denotes the set
of homogeneous elements with degree $j$. The \textit{ graded
dimension of }$\M$ is the power series
\[
q\dim \M:=\mathrel{\mathop{\sum }\limits_{j}}q^{j}\cdot rank(M_{j}).
\]
\noindent where $rank(M_j)=\dim_{\mathbb{Q}} (M_j \otimes
\mathbb{Q})$.

\end{definition}

\begin{remark}
$\M$ may have torsion but the \textit{graded dimension} will not
detect it.
\end{remark}


The direct sum and tensor product can be defined in the graded
category in an obvious way. The following proposition follows easily
from definition and the proof is omitted.

\begin{proposition} Let $\M\
$and $\N$ be a graded $\mathbb{Z}$-modules. Then $ \M\oplus \N$
and $\M\otimes \N$ are both graded $\mathbb{Z}$-modules with:
\begin{align*}&q\dim \left( \M\oplus \N\right) =q\dim \left( \M\right) +q\dim
\left( \N\right)\qquad\text{and}\\
&q\dim \left( \M\otimes \N\right) =q\dim \left( \M\right) \cdot
q\dim \left( \N\right) \end{align*}
\end{proposition}

\begin{example}
\label{M} Let $M$ be the graded free $\mathbb{Z}$-module with two
basis elements $1$ and $x$ \ whose degrees are 0 and 1
respectively. We have $M=\mathbb{Z} \oplus \mathbb{Z}x$ and $q
\dim M=1+q.$ This is the $\mathbb{Z}$-module we will use to
construct our chain complex. We have $q\dim \left( M^{\otimes
k}\right) =\left( 1+q\right) ^{k}.$ \\In order to avoid confusion,
we reserve the notation straight $M$ for this specific module and
use $\M$ for a generic module.
\end{example}


\begin{definition}Let $\left\{ \ell \right\} $ be the ``degree
shift'' operation on graded $\mathbb{Z}$-modules. That is, if
$\M={\oplus }_{j} M_{j}$ is a graded $\mathbb{Z}$-module where
$M_j$ denotes the set of elements of $\M$ of degree $j$, we set
$\M\left\{\ell \right\} _{j}:=M_{j-\ell }$ so that $\ q\dim
\M\left\{ \ell \right\} =q^{\ell }.q\dim \M.$ In other words, all
the degrees are increased by $\ell. $
\end{definition}
Therefore, $\mathbb{Z}\{ 1\}$ denotes $\mathbb{Z}$ with degree of
each element being 1. It is easy to check that $\M\otimes
\mathbb{Z}\{ 1\} \cong \M\{1\}$, the $ \mathbb{Z}$-module
isomorphic to $\M$ with degree of every homogeneous element raised
up by 1.

We now explain our construction. Let $G$ be a graph with $n$
ordered edges, and let $M$ be as in Example \ref{M}. For each
vertex $\varepsilon =(\varepsilon _{1},\varepsilon
_{2},....,\varepsilon _{n})$ of $\ $the cube $\left\{ 0,1\right\}
^{n},$ let $k_{\varepsilon }$ be the number of components of
$G_{\varepsilon }$. We assign a copy of $M$ to each connected
component and then take tensor product. This yields a graded free
$\mathbb{Z}$-module $ M_{\varepsilon}(G)=M^{\otimes k_{\varepsilon
}}$. Now, $q\dim M_{\varepsilon }(G)$ is the polynomial that
appears in the vertex $\varepsilon $ of the cube in 
\figref{figure State sum Chrom}.

\begin{figure}[ht!]\anchor{spaces}
\begin{center}
\cl{\scalebox{.65}{\includegraphics{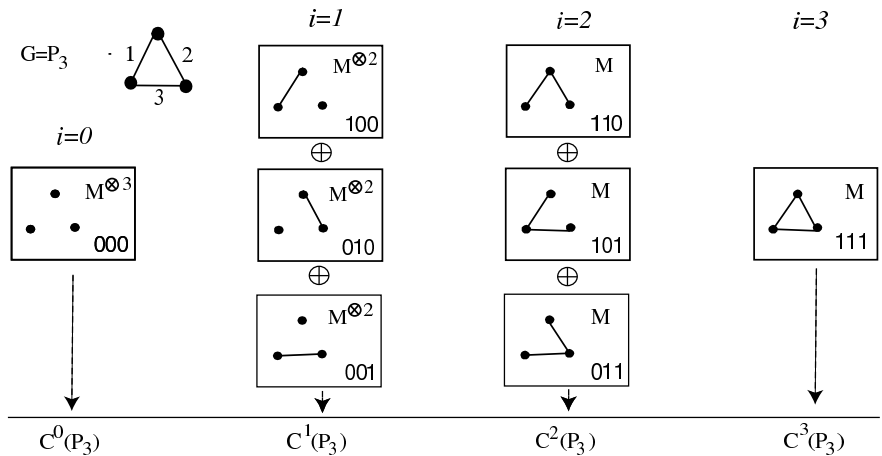}}}
\end{center}
\caption{The chain groups $C^{i}(G)$}
\label{spaces}
\end{figure}

To get the chain groups, we ``flatten'' the cube by taking direct
sums along the columns. A more precise definition is:

\begin{definition}
We set the $i$-th chain group $C^{i}(G)$ of the chain complex $\mathcal{C(}G%
\mathcal{)}$ to be the direct sum of all $\mathbb{Z}$-modules at height $i$, 
i.e.\ $C^{i}(G)=
\mathrel{\mathop{\oplus }\limits_{\left| \varepsilon \right| =i}} M_{\varepsilon }(G).$
\newline The grading is
given by the degree of the elements and  we can write the $i$-th
chain group as $C^{i}(G)= \oplus _{j\geq 0}C^{i,j}(G)$
 where $C^{i,j}(G)$ denotes the elements of degree $j$
of $\ C^{i}(G)$.
\end{definition}

\begin{figure}[ht!]\anchor{Basis elts C11(P3)}
\cl{\scalebox{.3}{\hglue-3cm\includegraphics{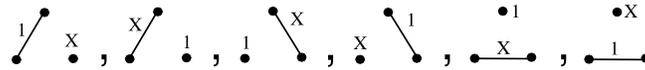}}}
\caption{Basis elements of the free $\mathbb{Z}$-module
$C^{1,1}(P_{3})$} \label{Basis elts C11(P3)}
\end{figure}

For example, the elements of degree one of $C^{1}(P_{3})$ are the
linear combinations with coefficients in $\mathbb{Z}$\ of the
 six elements shown in Figure \ref{Basis elts C11(P3)}.

These elements form a basis of the free $\mathbb{Z}$-module
$C^{1,1}(P_{3})$. This will lead to a second description of our
chain complex explained in section 2.3.

\paragraph{Graded chain complex, graded Euler characteristic}

\begin{definition}
Let $\M=\oplus_j M_j$ and $\N=\oplus_j N_j$ be graded
$\mathbb{Z}$-modules where $M_j$ (resp. $N_j$) denotes the set of
elements of $\M$ (resp. $\N$) of degree $j$.
\newline \noindent
A $\mathbb{Z}$-module map $\alpha: \M\rightarrow \N$ is said to be
\emph{graded with degree $d$} if $\alpha(M_j) \subseteq N_{j+d}$,
i.e.\ elements of degree $j$ are mapped to elements of degree
$j+d.$
\newline \noindent A \emph{graded chain complex} is a chain complex for which
the chain groups are graded $\mathbb{Z}$-modules and the
differentials are graded.
\end{definition}

\begin{definition}
The graded Euler characteristic $\chi_{q}(\mathcal{C})$ of a
graded chain complex $\mathcal{C}$ is the alternating sum of the
graded dimensions of its cohomology groups, i.e.\
$\chi_{q}(C)=\mathrel{\mathop{\sum }\limits_{0\leq i\leq
n}}(-1)^{i}\cdot q\dim (H^{i}).$
\end{definition}

The following was observed in \cite{BN02}. For convenience of the
reader, we include a proof here.

\begin{proposition}
\label{GEC}If the differential is degree preserving and all chain
groups are finite dimensional, the graded Euler characteristic is
also equal to the alternating sum of the graded dimensions of its
chain groups i.e. \[
\chi _{q}(C)=\mathrel{\mathop{\sum
}\limits_{0\leq i\leq n}}(-1)^{i}\cdot q\dim
(H^{i})=\mathrel{\mathop{\sum }\limits_{0\leq i\leq n}}
(-1)^{i}\cdot q\dim (C^{i}).
\]
\end{proposition}

\begin{proof}
The corresponding result for the non-graded case is well known.
That is, for a finite chain complex $\mathcal{C}=$ $0\rightarrow
C^{0}\rightarrow C^{1}\rightarrow ...\rightarrow C^{n}\rightarrow
0$ with cohomology groups $ H^{0},H^{1},...,H^{n},$ if all the
chain groups are finite dimensional then the Euler characteristic
$\chi (\mathcal{C})= \sum_{0\leq i\leq n}(-1)^{i} rank(H^{i})$ is
also equal to \\
$ \sum_{0\leq i\leq n}(-1)^{i} rank(C^{i}).$

Now, let $\mathcal{C}$ be a graded chain complex with a degree
preserving differential. With the above notations, decomposing
elements by degree yields $ C^{i}=\mathrel{\mathop{\oplus
}\limits_{j\geq 0}}C^{i,j}(G).$ Since the differential is degree
preserving, the restriction to elements of degree $j$, i.e.\
$0\rightarrow C^{0,j}\rightarrow C^{1,j}\rightarrow ...\rightarrow
C^{n,j}\rightarrow 0$ is a chain complex. The previous equation for
non-graded case tells us $ \mathrel{\mathop{\sum }\limits_{0\leq
i\leq n}}(-1)^{i}rank(H^{i,j})$ \\ $= \mathrel{\mathop{\sum
}\limits_{0\leq i\leq n}}(-1)^{i}rank(C^{i,j}).$ Now, multiply this
by $q^{j}$ and take the sum over all values of $j$ and the result
follows.
\end{proof}

\subsubsection{The differential}

\begin{figure}[ht!]
\cl{\scalebox{.7}{\includegraphics{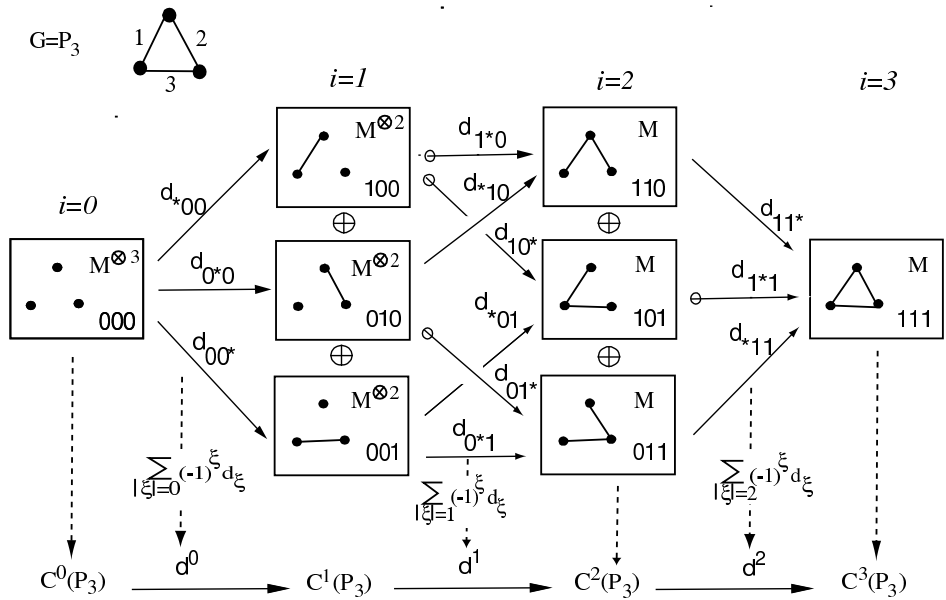}}}
\caption{The differentials} \label{The chain maps}
\end{figure}

\figref{The chain maps} shows what the maps look like. Details
are as follows.

We will first turn the edges of the cube $\left\{ 0,1\right\}^{n}$
into maps between the $\Z$-modules corresponding to its ends. We
call these maps \textit{per-edge maps}.


Recall that each vertex of the cube $\left\{ 0,1\right\} ^{n}$ is
labeled with some $\varepsilon=(\varepsilon_{1},....,\varepsilon
_{n})\in \left\{ 0,1\right\} ^{n}.$

$\blacktriangle$\qua\textsl{For which pairs of vertices are there
per-edge maps?}

There is a map between the $\Z$-modules corresponding to two
vertices if one of the markers $\varepsilon _{k}$ is changed from
$0$ to $1$ when you go from the first vertex to the second vertex
and all the other $\varepsilon _{k}$ are unchanged, and no map
otherwise.

Denote by $\varepsilon $ the label of the first vertex. If the
marker which is changed from $0$ to $1$ has index $k_{0}$ then the
map will be labelled $ d_{\varepsilon ^{\prime }}^{+}$  where
$\varepsilon ^{\prime }=(\varepsilon _{1}^{\prime
},....,\varepsilon _{n}^{\prime })$ with $\varepsilon _{k}^{\prime
}=$ $\varepsilon _{k}$ if $k\neq k_{0}$ and $ \varepsilon
_{k}^{\prime }=\ast $ if $k=k_{0}.$

For example, in the above diagram, the label $0\ast 1$ of the map
$d_{0\ast 1 \text{ }}^{+}$ means its domain is the $\Z$-module
corresponding to the vertex labelled $001$ and its target is the
$\Z$-modules corresponding to the vertex labelled $011.$ The
superscript $+$ indicates that we have not assigned signs to the
map yet.

$\blacktriangle $\qua\textsl{Definition of the per-edge maps}

Changing exactly one marker from $0$ to $1$ corresponds to adding an edge.

$\vartriangle $ If adding that edge doesn't change the number of
components, then the map is identity on $M^{\otimes k}.$

$\vartriangle $ If adding that edge decreases the number of
components by one, then we set $d_{\varepsilon ^{\prime
}}^{+}:M^{\otimes k}\rightarrow M^{\otimes k-1}$ to be identity on
the tensor factors corresponding to components that don't
participate and to be the $\mathbb{Z}$-linear map $m:$ $M\otimes
M$ $\rightarrow M$ given by $m(1\otimes 1)=1$, $m(1\otimes
x)=m(x\otimes 1)=x$ and $\ m(x\otimes x)=0$ on the two affected
components.

Note that identity and $m$ are degree preserving so $d_{\varepsilon
^{\prime }}^{+}$ inherits this property.

$\blacktriangle $\qua\textsl{``Flatten'' to get the differential}

The differential $d^{i}:C^{i}(G)\rightarrow C^{i+1}(G)$ of the chain
complex $\mathcal{C(}G\mathcal{)}$ is defined by
$d^{i}:=\mathrel{\mathop{\sum }\limits_{\left| \varepsilon \right|
=i}}$ $(-1)^{\varepsilon }d_{\varepsilon }^{+}$ where $\left|
\varepsilon \right| $ is the number of $ 1$'s in $\varepsilon $ and
$(-1)^{\varepsilon }$ is defined in the next paragraph.

$\blacktriangle$\qua \textsl{Assign a\ }$\pm 1$\textsl{ factor to each
per-edge map }$d_{\varepsilon }^{+}$

Equipped with the maps $d_{\varepsilon }^{+}$, the cube $\left\{
0,1\right\} ^{n}$ is commutative. This is because the multiplication
map $m$ is associative and commutative. To get the differential $d$
to satisfy $d\circ d=0,$ it is enough to assign a$\ \pm 1$ factor to
these maps in the following way: Assign $-1$ to the maps that have
an odd number of $1$'s before the star in
their label $\varepsilon ,$ and $1$ to the others. This is what was denoted $%
(-1)^{\varepsilon }$ in the definition of the differential.

In \figref{The chain maps}, we have indicated the maps for
which $(-1)^{\varepsilon }=-1$ by a little circle at the tail of
the arrow.

A straightforward calculation implies:

\begin{proposition}
This defines a differential, that is, $d^2 = 0$.
\end{proposition}

Now, we really have a chain complex $\mathcal{C(}G\mathcal{)}$ where
the chain groups and the differential are defined as above.
According to Proposition \ref{GEC}, we have:

\begin{theorem}
The Euler characteristic of this chain complex
$\mathcal{C(}G\mathcal{)}$ is equal to the chromatic polynomial of
the graph $G$ evaluated at $\lambda =1+q$.
\end{theorem}

\subsubsection{Independence of ordering of edges}

Let $G$ be a graph with edges labeled $1$ to $n$. For any
permutation $ \sigma $ of $\{1,..,n\},$ we define $G_{\sigma }$ to
be the same graph but with edges labeled as follows. The edge which
was labeled $k$ in $ G $ is labeled $\sigma (k)$ in $G_{\sigma }$.
In other words, $G$ is obtained from $G_{\sigma }$ by permuting the
labels of the edges of $G$ according to $ \sigma .$

\begin{theorem}
\label{ordering} The chain complexes $\mathcal{C(}G\mathcal{)}$
and $\mathcal{C(}G\mathcal{ _{\sigma })}$ are isomorphic and
therefore, the cohomology groups are isomorphic.  In other words,
the cohomology groups are independent of the ordering of the edges
so they are well defined graph invariants.
\end{theorem}

{\bf Proof.} Since the group of permutations on $n$ elements is generated by the permutations of the form
$(k,k+1)$ it is enough to prove the result when $\sigma =(k,k+1).$

We will define an isomorphism $f$ such that the following diagram
commutes:

$\hspace{1cm}\mathcal{C(}G\mathcal{)}$ : \ $C^{0}(G)\overset{d^{0}}{\
\rightarrow } C^{1}(G)\overset{d^{1}}{\rightarrow
}...\overset{d^{i-1}}{\rightarrow }
C^{i}(G)\overset{d^{i}}{\rightarrow
}...\overset{d^{n-1}}{\rightarrow } C^{n}(G)$

\hspace{1cm}\ \ \ \ \ \ \ \ \ \ \ \ $\downarrow f$ \ $\ \ \ \ \ \ \ \ \downarrow f$ \ $\
\ \ \ \ \ \ \ \ \ \ \ \ \ \ \ \ \downarrow f$ \ $\ \ \ \ \ \ $\ \ ...$\ \ \
\ \ \downarrow f$ \ \

$\hspace{1cm}\mathcal{C(}G\mathcal{_{\sigma })}$: $C^{0}(G\mathcal{_{\sigma
}})\overset{ d^{0}}{\rightarrow }C^{1}(G\mathcal{_{\sigma
}})\overset{d^{1}}{\rightarrow } ...\overset{d^{i-1}}{\rightarrow
}C^{i}(G\mathcal{_{\sigma }})\overset{d^{i} }{\rightarrow
}...\overset{d^{n-1}}{\rightarrow }C^{n}(G\mathcal{_{\sigma }} ) $

It is enough to define $f$ restricted on each submodule
$M_{\varepsilon}(G)$ (where $|\varepsilon |=i$)  since $C^{i}(G)$
is the direct sum of these submodules.

\begin{figure}[ht!]\anchor{permutLabels}
\begin{center}
\cl{\scalebox{.75}{\includegraphics{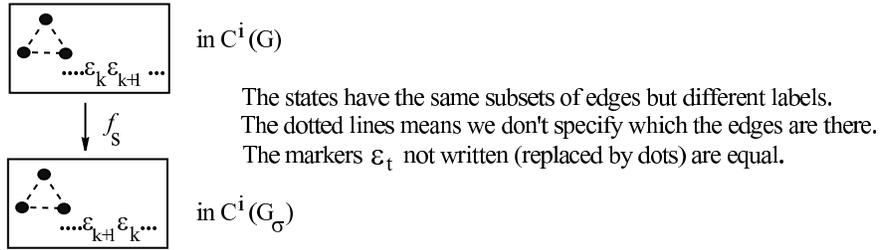}}}
\end{center}
\caption{Impact of re-ordering of the edges}
\label{permutLabels}
\end{figure}

For any subset $s$ of $E$ with $i$ edges, there is a state in $C^{i}(G)$ and
one in $C^{i}(G\mathcal{_{\sigma }})$ that correspond exactly to those
edges. Let $\varepsilon =(\varepsilon _{1},....,\varepsilon _{n})$ stand for
$\varepsilon _{s}(G),$ the label of $s$ in $G.$ The situation is illustrated
in \figref{permutLabels}.

Let $f_{s}$ be the map between these two states that is equal to
$-id$ if $ \varepsilon _{k}=\varepsilon _{k+1}=1$ and equal to $id$
otherwise.

Let $f:C^{i}(G)\rightarrow C^{i}(G\mathcal{_{\sigma }})$ be defined
by $f= \mathrel{\mathop{\oplus }\limits_{\left| s\right| =i}}f_{s}.$

$f$ is obviously\ an isomorphism. The fact that the diagram
commutes can be checked by looking at the four cases $(\varepsilon
_{k},\varepsilon _{k+1})=(0,$ $0)$, $(\varepsilon _{k},\varepsilon
_{k+1})=(1,$ $0)$, $ (\varepsilon _{k},\varepsilon _{k+1})=(0,$
$1)$ and $(\varepsilon _{k},\varepsilon _{k+1})=(1,$ $1)$.


\subsubsection{A Poincar\'{e} polynomial}
We define a two-variable polynomial $R_{G}(t,q)$ by $R_{G}(t,q)\,{=}\!\!\!\!
\mathrel{\mathop{\sum }\limits_{ 0\leq i\leq n}}\!\!\!\!t^{i}\cdot q\dim
\left( H^i (G)\right)$. In the following proposition, (a) follows
immediately from Theorem \ref{ordering} and (b) follows from our
construction.

\begin{proposition}
{\rm(a)}\qua The polynomial $R_{G}(t,q)$ depends only on the graph.

{\rm(b)}\qua The chromatic polynomial is a specialization of
$R_{G}(t,q)$ at $t=-1.$
\end{proposition}

This polynomial is a convenient way to store the information about
the free part of the cohomology groups and is, by construction,
enough to recover the chromatic polynomial.

\subsection{Another description: The enhanced state construction}

These cohomology groups have another description that is similar
to Viro's description for the Khovanov cohomology for knots
\cite{V02}. We explain the details below.

Let $\{1,x\}$ be a set of \textit{colors}, and $\ast $ be the
product defined by
\begin{equation*}
1\ast 1=1,1\ast x=x\ast 1=x\text{ and \ }x\ast x=0
\end{equation*}
Let $G=(V,E)$ be a graph with an ordering on its edges. An
\textit{enhanced state} of $G$ is $S=(s,c)$, where $s\subseteq E$
and $c$ is an assignment of $1$ or $x$ to each connected component
of the spanning subgraph $[G:s]$. For each enhanced state $S$,
define
\begin{equation*}
i(S)=\#\text{ of edges in }s,\text{ and }j(S)=\#\text{ of }x\text{
in }c.
\end{equation*}
Note that $i(S)$ depends only on the underlying state $s,$ not on
the color assignment that makes it an enhanced state, so we may
write it as $i(s).$

Let $C^{i,j}(G):=Span \langle S|S$ is an enhanced state of $G$ with
$i(S)=i,j(S)=j \rangle$, where the span is taken over $\mathbb{Z}$.

We define the differential
\begin{equation*}
d:C^{i.j}(G)\rightarrow C^{i+1,j}(G)
\end{equation*}
as follows. For each enhanced state $S=(s,c)$ in $C^{i,j}(G)$, define $%
d(S)\in C^{i+1,j}(G)$ by
\begin{equation*}
d(S)=\sum_{e\in E(G)-s}(-1)^{n(e)}S_{e}
\end{equation*}
\noindent where $n(e)$ is the number of edges in $s$ that are ordered before
$e$, $S_{e}$ is an enhanced state or 0 defined as follows. Let $s_{e}=s\cup
\{e\}$. Let $E_{1},\cdots ,E_{k}$ be the components of $[G:s]$. If $e$
connects some $E_{i}$, say $E_{1}$, to itself, then the components of $%
[G:(s\cup \{e\})]$ are $E_{1}\cup \{e\},E_{2},\cdots ,E_{k}$. We define $%
c_{e}(E_{1}\cup \{e\})=c(E_{1}),c_{e}(E_{2})=c(E_{2}),\cdots
,c_{e}(E_{k})=c(E_{k})$, and $S_{e}$ is the enhanced state $(s_{e},c_{e})$.
If $e$ connects some $E_{i}$ to $E_{j}$, say $E_{1}$ to $E_{2}$, then the
components of $[G:s_{e}]$ are $E_{1}\cup E_{2}\cup \{e\},E_{3},\cdots ,E_{k}$%
. We define $c_{e}(E_{1}\cup E_{2}\cup \{e\})=c(E_{1})\ast
c(E_{2}),c_{e}(E_{3})=c(E_{3}),\cdots ,c_{e}(E_{k})=c(E_{k})$. Note that if $%
c(E_{1})=c(E_{2})=x, c_{e}(E_{1}\cup E_{2}\cup \{e\})= x\ast x =0$, and therefore $%
c_{e} $ is not considered as a coloring. In this case, we let $S_{e}=0$. In
all other cases, $c_{e}$ is a coloring and we let $S_{e}$ be the enhanced
state $(s_{e},c_{e})$.

One may find it helpful to think of $d$ as the operation that adds
each edge not in $s$, adjusts the coloring using $\ast $, and then
sums up the enhanced states using appropriate signs. In the case
when an illegal color of 0 appears, due to the product $x\ast
x=0$, the contribution from that edge is counted as 0.

\subsection{Equivalence of the two constructions}

At first sight, the two constructions look different because the
cubic complex construction yields only one chain complex whereas the
enhanced states construction gives rise to a sequence of chain
complexes, one for each degree $j$. This can be easily solved by
splitting the chain complex of the cubic complex construction into a
sequence of chain complexes, one for each degree $j$. More
precisely, let $\mathcal{C}$ = $0\rightarrow C^{0}\rightarrow
C^{1}\rightarrow ...\rightarrow C^{n}\rightarrow 0$ be a graded
chain complex with a degree preserving differential. Decomposing
elements of each chain group by degree yields
$C^{i}=\mathrel{\mathop{\oplus }\limits_{j\geq 0}}C^{i,j}.$ Since
the differential is degree preserving, the restriction to elements
of degree $j$, i.e.\ $0\rightarrow C^{0,j}\rightarrow
C^{1,j}\rightarrow ...\rightarrow C^{n,j}\rightarrow 0$ is a chain
complex denoted  by $\mathcal{C}^{j}$.

We are now ready to see that for a fixed $j$, the chain complexes obtained
via the two construction are isomorphic. For this paragraph, denote the one
obtained via the cubic complex construction by $\mathcal{C}^{j}$ and the one
obtained via the enhanced state construction by $\widetilde{\mathcal{C}}%
^{j}. $

Both chain complexes have free chain groups so it is enough to
define the chain map on basis elements. We will associate each
enhanced state $S=(s,c)$ of $\widetilde{C}^{i,j}(G)$ to an unique
basis element in $C^{i,j}(G)$ and show that this defines an
isomorphism of chain complexes. First, $s\subseteq E(G)$ naturally
corresponds to the vertex $\varepsilon =(\varepsilon _{1},\cdots
,\varepsilon _{n})$ of the cube, where $\varepsilon _{k}=1$ if $
e_{k}\in s$ and $\varepsilon _{k}=0$ otherwise. The corresponding
$\mathbb{Z} $-module $M_{\varepsilon }(G)$ is obtained by assigning
a copy of $M$ to each connected component of$\ [G,s]$ and then
taking tensor product. The color $c$ naturally corresponds to the
basis element $x_{1}\otimes \cdots \otimes x_{k}$ where $x_{\ell }$
is the color associated to the $\ell$-th component of $[G:s]$ .

It is not difficult to see that this defines an isomorphism on the chain
group that commutes with the differentials. Therefore, the two complexes are
isomorphic.

\section{Properties}

In this section, we demonstrate some properties of our cohomology
theory, as well as some computational examples.

\subsection{An exact sequence}

The chromatic polynomial satisfies a well-known deletion-contraction
rule: $$ P(G,\lambda )=P(G-e,\lambda )-P(G/e,\lambda ).$$
 Here we show that
our cohomology groups satisfy a naturally constructed long exact
sequence involving $G,G-e$ , and $G/e$. Furthermore, by taking the
Euler-characteristic of the long exact sequence, we recover the
deletion-contraction rule. Thus our long exact sequence can be
considered as a categorification of the deletion-contraction rule.

 We explain the exact sequence in terms of the enhanced
state sum approach. Let $G$ be a graph and $e$ be an edge of $G$. We
order the edges of $G$ so that $e$ is the last edge. This induces
natural orderings on $G/e$ and on $G-e$ by deleting $e$ from the
list. We define homomorphisms $\alpha
_{ij}:C^{i-1,j}(G/e)\rightarrow C^{i,j}(G)$ and $\beta
_{ij}:C^{i,j}(G)\rightarrow C^{i,j}(G-e)$. These two maps will be
abbreviated by $\alpha $ and $\beta $ from now on. Let $v_{e}$ and
$w_{e}$ be the two vertices in $G$ connected by $e$. Intuitively,
$\alpha $ expands $v_{e}=w_e$ by adding $e$, and $\beta $ is the
projection maps. We explain more details below.

 First, given an enhanced state $S=(s,c)$ of $G/e$, let $
\widetilde{s}=s\cup \{e\}$.  The number of components of $[G/e:s]$
and $[G:\widetilde{s}]$ are the same. In  fact, the components of
$[G/e:s]$ and the components of $[G:\widetilde{s}]$ are the same
except the one containing $v_{e}$ where $ v_{e}$ in $G/e$ is
replaced by $e$ in $G$. Thus, $c$ automatically yields a coloring of
components of $[G:\widetilde{s}]$, which we denote by $
\widetilde{c}$. Let $\alpha (S)=(\widetilde{s},\widetilde{c})$. It
is an enhanced state in $C^{i,j}(G)$. Extend $\alpha $ linearly and
we obtain a homomorphism $\alpha :C^{i-1,j}(G/e)\rightarrow
C^{i,j}(G)$.

 Next, we define the map $\beta :C^{i,j}(G)\rightarrow
C^{i,j}(G-e)$. Let $S=(s,c)$ be an enhanced state of $G$. If $e\not\in s$, $S
$ is automatically an enhanced state of $G-e$ and we define $\beta (S)=S$.
If $e\in s$, we define $\beta (S)=0$. Again, we extend $\beta $ linearly to
obtain the map $\beta :C^{i,j}(G)\rightarrow C^{i,j}(G-e)$.

One can sum up over $j$, and denote the maps by $\alpha_i:
C^{i-1}(G/e)\rightarrow C^i(G)$ and $\beta_i: C^i(G)\rightarrow
C^i(G-e)$. Again, they will be abbreviated by $\alpha $ and $\beta$.
Both are degree preserving maps since the index $j$ is preserved.

\begin{lemma}
$\alpha $ and $\beta$ are chain maps such that $0\rightarrow
C^{i-1,j}(G/e)\overset{\alpha}{\rightarrow} C^{i,j}(G)
\overset{\beta}{\rightarrow} C^{i,j}(G-e) \rightarrow 0$ is a short
exact sequence.
\end{lemma}

\begin{proof}
First we show that $\alpha $ is a chain map. That is,
\begin{equation*}
C^{i-1,j}(G/e)\overset{\alpha }{\rightarrow }C^{i,j}(G)
\end{equation*}
\begin{equation*}
\downarrow d_{G/e}\ \ \ \ \ \ \ \ \ \ \ \downarrow d_{G}
\end{equation*}
\begin{equation*}
C^{i,j}(G/e)\overset{\alpha }{\rightarrow }C^{i+1,j}(G)
\end{equation*}
commutes. Let $(s,c)$ be an enhanced state of $G/e$, we have
\begin{equation*}
d_{G}\circ \alpha ((s,c))=d_{G}(s\cup
\{e\},\widetilde{c})=\sum_{e_{i}\in E(G)-(s\cup
\{e\})}(-1)^{n_{G}(e_{i})}\left( s\cup \{e,e_{i}\},(\widetilde{c}
)_{e_{i}}\right)
\end{equation*}
where $n_{G}(e_{i})$ is the number of edges in $s\cup \{ e \}$ that
are ordered before $e_{i}$ in $G$, and $(\widetilde{c} )_{e_{i}}$ is
the natural coloring inherited from $\widetilde{c}$ when adding the
edge $e_i$ (see Section 2.3 for the description).

We also have
 $\alpha \circ d_{G/e}((s,c))=\alpha \left(
\sum_{e_{i}\in E(G/e)-s}(-1)^{n_{G/e}(e_{i})}(s\cup \{ e_{i}
\},c_{e_{i}})\right)$ \\
 $=\sum_{e_{i}\in
E(G/e)-s}(-1)^{n_{G/e}(e_{i})}(s\cup
\{e_{i},e\},\widetilde{(c_{e_{i}})})$
 where $n_{G/e}(e_{i})$ is
the number of edges in $s$ that are ordered before $e_{i}$ in $G/e$.

The two summations contain the same list of $e_{i}$'s since
$E(G)-(s\cup \{e\})=E(G/e)-s$. It is also easy to see that
$(\widetilde{c})_{e_{i}}= \widetilde{(c_{e_{i}})}$. Finally,
$n_{G}(e_{i})=n_{G/e}(e_{i})$ since $e$ is ordered last. It follows
that $d_{G}\circ \alpha =\alpha \circ d_{G/e}$ and therefore $\alpha
$ is a chain map.

 Next, we show that $\beta $ is a chain map by proving the
commutativity of
\begin{equation*}
C^{i,j}(G)\overset{\beta}{\rightarrow } C^{i,j}(G-e)
\end{equation*}
\begin{equation*}
\downarrow d_G \ \ \ \ \ \ \ \ \ \ \ \downarrow d_{G-e}
\end{equation*}
\begin{equation*}
C^{i+1,j}(G)\overset{\beta}{\rightarrow } C^{i+1,j}(G-e)
\end{equation*}
Let $S=(s,c)$ be an enhanced state of $G$.

If $e\in s$, we have $\beta (S)=0$ and thus $d_{G-e}\circ \beta
(S)=0$. We also have $d_{G}(S)=\mathrel{\mathop{\sum
}\limits_{e_{i}\in E(G)-s}}(-1)^{n_{G}(e_{i})}(s\cup
\{e_{i}\},c_{e_{i}})$. Since $e\in s\cup \{e_{i}\}, $ $\beta \circ
d_{G}(S)=0$.

If $e\not\in s$, we have $d_{G-e}\circ \beta (S)=d_{G-e}(S)=
\mathrel{\mathop{\sum }\limits_{e_{i}\in
E(G-e)-s}}(-1)^{n_{G-e}(e_{i})}(s\cup \{e_{i}\},c_{e_{i}}) $. We
also have $d_{G}(S)=\mathrel{\mathop{ \sum }\limits_{e_{i}\in
E(G)-s}}(-1)^{n_{G}(e_{i})}(s\cup
\{e_{i}\},c_{e_{i}})=S_{1}+S_{2}$, where $S_{1}=%
\mathrel{\mathop{\sum
}\limits_{e_{i}\in E(G)-(s\cup \{e\})}}(-1)^{n_{G}(e_{i})}(s \cup
\{e_{i}\},c_{e_{i}})$ corresponds to the terms with $e_{i}\neq e$, and $%
S_{2}=(-1)^{n_{G}(e_{i})}(s\cup \{e\},c_{e})$ corresponds to the
term $ e_{i}=e$. By our definition of $\beta $, $\beta
(S_{1})=S_{1},\beta (S_{2})=0 $. Finally,
$n_{G}(e_{i})=n_{G-e}(e_{i})$ since $e$ is ordered last, and it
follows that $d_{G-e}\circ \beta (S)=\beta \circ d_{G}(S)$ as well
in this case.

Finally, we prove the exactness. Each nonzero element in $
C^{i-1,j}(G/e)$ can be written as $x=\sum n_{k}(s_{k},c_{k})$
where $ n_{k}\neq 0$ and $(s_{k},c_{k})$'s are distinct enhanced
states of $ G/e$. It is not hard to see that
$(\tilde{s_{k}},\tilde{c_{k}})$'s are distinct enhanced states of
$G$. Thus $\alpha (x)=\sum n_{k}(\tilde{ s_{k}},\tilde{c_{k}})\neq
0$ in $C^{i-1,j}(G)$. Hence $\ker\alpha =0$. Next,
$\mbox{Im}\alpha =\ker \beta$, since both are
 Span $\{(s,c)\:|\:(s,c)$ is an enhanced state of $ G$ and $e\in
s\}$. Last, $\beta $ is a projection map that maps onto
$C^{i,j}(G-e)$.
\end{proof}

The Zig-Zag lemma in homological algebra implies :

\begin{theorem}
Given a graph $G$ and an edge $e$ of $G$, for each $j$ there is a long exact
sequence

$0\rightarrow H^{0,j}(G)\overset{\beta ^{\ast }}{\rightarrow
}H^{0,j}(G-e)
\overset{\gamma ^{\ast }}{\rightarrow }H^{0,j}(G/e)\overset{\alpha ^{\ast }}{%
\rightarrow }H^{1,j}(G)\overset{\beta ^{\ast }}{\rightarrow }H^{1,j}(G-e)%
\overset{\gamma ^{\ast }}{\rightarrow }H^{1,j}(G/e)\rightarrow \ldots
\rightarrow H^{i,j}(G)\overset{\beta ^{\ast }}{\rightarrow }H^{i,j}(G-e)%
\overset{\gamma ^{\ast }}{\rightarrow }H^{i,j}(G/e)\overset{\alpha ^{\ast }}{%
\rightarrow }H^{i+1,j}(G)\rightarrow \ldots \ $

If we sum over $j$, we have a degree preserving long exact sequence:

$0\rightarrow H^{0}(G)\overset{\beta ^{\ast }}{\rightarrow
}H^{0}(G-e) \overset{\gamma ^{\ast }}{\rightarrow
}H^{0}(G/e)\overset{\alpha ^{\ast }}{ \rightarrow
}H^{1}(G)\overset{\beta ^{\ast }}{\rightarrow }H^{1}(G-e)
\overset{\gamma ^{\ast }}{\rightarrow }H^{1}(G/e)\rightarrow \ldots
\rightarrow H^{i}(G)\overset{\beta ^{\ast }}{\rightarrow
}H^{i}(G-e)\overset{ \gamma ^{\ast }}{\rightarrow
}H^{i}(G/e)\overset{\alpha ^{\ast }}{ \rightarrow
}H^{i+1}(G)\rightarrow \ldots $
\end{theorem}

\begin{remark}
\label{map description} It is useful to understand how the maps
$\alpha ^{\ast },\beta ^{\ast },\gamma ^{\ast }$ act in an
intuitive way. The descriptions for $\alpha $ and $\beta $ follows
directly from our construction: $\alpha ^{\ast }$ expands the edge
$e$, $\beta ^{\ast }$ is the projection map. The description for
$\gamma ^{\ast }$ follows from the standard diagram chasing
argument in the zig-zag lemma and the result is as follows. For
each cycle $z$ in $C^{i,j}(G-e)$ represented by the chain $\sum
n_{k}(s_{k},c_{k}) $, $\gamma ^{\ast }(z)$ is represented by the
chain $ (-1)^{i}\sum n_{k}(s_{k}\cup \{e\}/e,(c_{k})_{e})$, where
$s_{k}\cup \{e\}/e$ is the subset of $E(G/e)$ obtained by adding
$e$ to $s_{k}$ and then contracting $e$ to $v_{e}$, $(c_{k})_{e}$
is the coloring defined in Section 2.3. We leave it to the reader
to verify the result.
\end{remark}

\subsection{Graphs with loops or multiple edges}

We prove two propositions.

\begin{proposition}
\label{graphs with loops} If the graph has a loop then all the
cohomology group are trivial.
\end{proposition}

\begin{proof}
Let $G$ be a graph with a loop $\ell$. The exact sequence for $(G,
\ell)$ is

$0\rightarrow H^{0}(G)\rightarrow H^{0}(G-\ell) \overset{\gamma
^{\ast }}{\rightarrow }H^{0}(G/\ell)\rightarrow H^{1}(G)\rightarrow
H^{1}(G-\ell) \overset{\gamma ^{\ast }}{\rightarrow
}H^{1}(G/\ell)\rightarrow \ldots \rightarrow H^{i}(G)\rightarrow
H^{i}(G-\ell)\overset{ \gamma ^{\ast }}{\rightarrow
}H^{i}(G/\ell)\rightarrow H^{i+1}(G)\rightarrow \ldots $

Using our description of the snake map $\gamma^{\ast}$ in Remark
\ref{map description}, we get that the map
$H^{i}(G-\ell)\overset{\gamma ^{\ast }}{\rightarrow } H^{i}(G/\ell)$
is $(-1)^i id$. Therefore, $H^{i}(G)=0$ for all $i$.
\end{proof}

\begin{proposition}
\label{graphs with multiple edges} The cohomology group are
unchanged if you replace all the multiple edges of a graph by
single edges.
\end{proposition}

\begin{proof}
Assume that in some graph $G$ the edges $e_{1}$\ and $e_{2}$
connect the same vertices. In $G/e_{2}$, $e_{1}$ becomes a loop so
as observed earlier, $ H^{i}(G/e_{2})=0$ for all $i$. It follows
from the long exact sequence that $ H^{i}(G-e_{2})$ and $H^{i}(G)$
are isomorphic groups. One can repeat the process until there is
at most one edge connecting two given vertices without changing
the cohomology groups.
\end{proof}

\subsection{Disjoint union of two graphs}

Let $G_{1}$\ and $G_{2}$\ be two graphs and consider their disjoint
union $ G_{1}\sqcup G_{2}.$\ On the chain complex level, we have
$\mathcal{C} (G_{1}\sqcup G_{2})=\mathcal{C}(G_{1})\otimes
\mathcal{C}(G_{2})$.

\begin{theorem}
\label{disjoint union}For each $i\in \mathbb{N}$, we have :
\begin{equation*}
H^{i}(G_{1}\sqcup G_{2})\cong \left[ \mathrel{\mathop{\oplus
}\limits_{p+q=i}} H^{p}(G_{1})\otimes H^{q}(G_{2})\right] \oplus
\left[ \mathrel{\mathop{ \oplus }\limits_{p+q=i+1}}
H^{p}(G_{1})*H^{q}(G_{2}) \right]
\end{equation*}
where * denotes the torsion product of two abelian groups.

\noindent If we decompose the groups by degree,
 we have the following for all $i, j\in \mathbb{N}$:
\begin{equation*}
H^{i,j}(G_{1}\sqcup G_{2})\cong \left[\!\! \mathrel{\mathop{\oplus
}\limits_{
\begin{array}{c}
{\scriptscriptstyle p+q=i} \\
{\scriptscriptstyle s+t=j}
\end{array}
}}\!\!H^{p,s}(G_{1})\otimes H^{q,t}(G_{2})\right]\oplus \left[\!\!
\mathrel{\mathop{\oplus }\limits_{
\begin{array}{c}
{\scriptscriptstyle p+q=i+1} \\
{\scriptscriptstyle s+t=j}
\end{array}\!\!
}} H^{p,s}(G_{1})* H^{q,t}(G_{2})\right]
\end{equation*}
\end{theorem}

\begin{proof}
This is a corollary of K\"{u}nneth theorem, since the chains
complexes $ \mathcal{C}(G_{1})$\ and $\mathcal{C}(G_{2})$\ are
free. See \cite{M84} for details about the K\"{u}nneth
theorem.
\end{proof}

\begin{corollary}
The Poincar\'{e} polynomials are multiplicative under disjoint union
i.e.\ $ R_{G_{1}\sqcup G_{2}}(t,q)=R_{G_{1}}(t,q)\cdot
R_{G_{2}}(t,q)$
\end{corollary}

Theorem \ref{disjoint union} also implies

\begin{example}
Disjoint union with the one vertex graph: 
$$H^{i}(G\sqcup
\bullet )\cong H^{i}(G)\otimes (\mathbb{Z}\oplus \mathbb{Z}x).$$
\end{example}

\subsection{Adding or contracting a pendant edge}

An edge in a graph is called a \textit{pendant} edge if the degree
of one of its end points is one. Let $G$ be a graph and $e$ be a
pendant edge of $G$. Let $G/e$ be the graph obtained by
contracting $e$ to a point. We will study the relation between the
cohomology groups of $G$ and $G/e$.

Recall that, for a given graded $\mathbb{Z}$-module $M$, $M\{1\}$
denotes the $\mathbb{Z}$-module isomorphic to $M$ with degree for
each homogeneous element being shifted up by 1. We have

\begin{theorem}
\label{pendant edge} Let $e$ be a pendant edge in a graph $G$.
\\ For each $i$, $H^{i}(G)\cong H^{i}(G/e)\{1\}$.
\end{theorem}

\begin{proof}
Consider the operations of contracting and deleting $e$ in $G$. For
convenience, denote the graph $G/e$ by $G_1$. We have $G/e=G_1,$ and $%
G-e=G_1\sqcup \{v\}$, where $v$ is the end point of $e$ with $\deg v
= 1$. The exact sequence on $(G,e)$ is
$\hspace{-1cm}
0\rightarrow H^{0}(G)\rightarrow H^{0}(G_1\sqcup \{v\})\rightarrow
H^{0}(G_1)\rightarrow \cdots \rightarrow H^{i}(G)\rightarrow
H^{i}(G_1 \sqcup \{v\})\rightarrow H^{i}(G_1)\rightarrow \cdots
$

Thus we need to understand the map
\begin{equation*}
H^{i}(G_1\sqcup \{v\})\overset{\gamma^*}{\rightarrow }H^{i}(G_1)
\end{equation*}
By Theorem $\ref{disjoint union}$,
\begin{equation*}
H^{i}(G_{1}\sqcup \{v\})\cong H^{i}(G_{1})\otimes \left[
\mathbb{Z}\oplus \mathbb{Z}\{1\}\right] \cong H^{i}(G_{1})\oplus
\left[ H^{i}(G_{1})\otimes \mathbb{Z}\{1\}\right]
\end{equation*}
by a natural isomorphism $h_*$, which is induced by the isomorphism
$h$ described as follows. Each enhanced state $S$ in
$C^{i}(G_{1}\sqcup \{v\})$ either assigns the color $1$ to $v$, or
the color $x$ to $v$. If it assigns $ 1$ to $v$, $h$ sends $S$ to\
$(-1)^{\left| s\right| }(S_{1},0)$ where $S_{1}$ is the
``restriction'' of $S$ to $G_{1}$ and $\left| s\right| $ is the
number of edges of the underlying set $s$ of $S=(s,c)$. If it
assigns $x$ to $v$, $h$ sends $S$ to $(-1)^{\left| s\right|
}(0,S_{1}\otimes x)$. This extends to a degree preserving
isomorphism on chain groups and induces the isomorphism $h^{\ast }$
on cohomology groups.

We therefore will identify $H^{i}(G_{1}\sqcup \{v\})$ with $
H^{i}(G_{1})\oplus \left[ H^{i}(G_{1})\otimes
\mathbb{Z}\{1\}\right] $.

\medskip
\textbf{\noindent Claim}\qua $\gamma ^{\ast }\circ (h^{\ast
})^{-1}:H^{i}(G_{1})\oplus \left[ H^{i}(G_{1})\otimes
\mathbb{Z}\{1\}\right] \rightarrow H^{i}(G_{1})$ satisfies \\
$\gamma ^{\ast }\circ (h^{\ast })^{-1}(x,0)=x\text{ for all }x\in
H^{i}(G_{1}).$

\medskip
\textbf{Proof of Claim}\qua
Let $x$ be in $H^{i}(G_{1})$. $x$ is the equivalence class of a sum
of terms of the form $(s,c)$ in $G_{1}$.  Under the map $(h^{\ast
})^{-1}$, each of these terms is ``extended'' to be an element in
$C^{i}(G_{1}\sqcup \{v\})$ by adding $v$ to $\left[ G_{1}:s\right] $
and assigning it the color $1$. The final result is then multiplied
by $(-1)^{i}$.

The map $\gamma ^{\ast }$ is described in Remark \ref{map
description}. Let $ y$ be in $H^{i}(G_{1}\sqcup \{v\})$. $y$ is the
equivalence class of a sum of terms of the form $(s,c)$ in
$H^{i}(G_{1}\sqcup \{v\})$. Basically, for each each $(s,c)$,
$\gamma ^{\ast }$ adds the edge $e$, adjusts the colorings, then
contracts $e$ to a point and multiplies the result by $ (-1)^{i}$.
Hence applying $(h^{\ast })^{-1}$ then $\gamma ^{\ast }$ yields the
original graph $G_{1}$. The color for each state in $x$ remains the
same since $v$ is colored by $1$ and multiplication by $1$ is the
identity map. This proves that $\gamma ^{\ast }\circ (h^{\ast
})^{-1}(x,0)=x$.

The claim implies that $\gamma ^{\ast }$ is onto for each $i$. Thus
the above long exact sequence becomes a collection of short exact
sequences.
\begin{equation*}
0\rightarrow H^{i}(G)\rightarrow H^{i}(G_{1}\sqcup \{v\})\rightarrow
H^{i}(G_{1})\rightarrow 0
\end{equation*}
After passing to the isomorphism $h_*$, the exact sequence becomes
\begin{equation*}
0\rightarrow H^{i}(G)\rightarrow H^{i}(G_{1})\oplus \left[
H^{i}(G_{1})\otimes \mathbb{Z}\{1\}\right] \rightarrow
H^{i}(G_{1})\rightarrow 0
\end{equation*}
The next lemma implies that $H^{i}(G)\cong H^{i}(G)\otimes
\mathbb{Z}\{1\}$.
\end{proof}

\begin{lemma}
Let $A$ and $B$ be graded abelian groups, and $p:A\oplus B\rightarrow A$ be
a degree preserving projection with $p(a,0)=a$ for all $a\in A$. Then $\ker
p\cong B$ via a degree preserving isomorphism.
\end{lemma}

\begin{proof}
For each $b\in B$, let $a_{b}=p(0,b)\in A$. Then $p(-a_{b},b)=0$ and
therefore $(-a_{b},b)\in \ker p$. Define $f(b)=(-a_{b},b)$. It is a
standard exercise to verify that $f$ is a degree preserving
isomorphism from $B$ to $ \ker p$.
\end{proof}

\subsection{Trees, circuit graphs}

We describe the cohomology groups for several classes of graphs.

\begin{example}
\label{1 vertex and no edges}Let $N_{1}$ be the graph with 1
vertex and no edge. Its chromatic polynomial is $P_{N_{1}}=\lambda
=1+q$. The only enhanced states of $N_{1}$ are $(\emptyset ,1)$
and $(\emptyset ,x)$, which generate $C^{0}(N_{1})$. It follows
that $H^{0}(N_{1})\cong \mathbb{Z}\oplus \mathbb{Z}\{ 1\}$, and $
H^{i}(N_{1})=0$ for all $i\neq 0.$
\end{example}

\begin{example}
\label{n vertices and no edges} More generally, the graph with $v$
vertices and no edges is called the \textit{null graph} of order
$v$ and denoted by $N_{v}.$ A similar argument implies
\begin{equation*}
H^{i}(N_{v})\cong \left\{
\begin{array}{ll}
(\mathbb{Z}\oplus \mathbb{Z}\{1\})^{\otimes v} & \text{if }i=0 \\
0 & \text{if }i\neq 0
\end{array}
\right.
\end{equation*}
This also follows from the K\"{u}nneth type formula in
Theorem $ \left( \ref{disjoint union}\right) $.
\end{example}


\begin{example}
\label{tree with n edges} Let $G=T_{n}$, a tree with $n$ edges. We
can obtain $G$ by starting from a one point graph, and then adding
pendant edges successively. Thus Theorem \ref{pendant edge} and
Example \ref{1 vertex and no edges} imply
\begin{equation*}
H^{0}(T_{n})\cong (\mathbb{Z}\oplus \mathbb{Z}\{1\})\{n\}\cong
\mathbb{Z} \{n\}\oplus \mathbb{Z}\{n+1\},\text{
}H^{i}(T_{n})=0\text{ for }i\neq 0.
\end{equation*}
\end{example}

A basis for $H^{0}(T_{n})$ can be described as follows. Let $
V(T_n)=\{v_{0},v_{1},\cdots ,v_{n}\}$ be the set of vertices of
$T_n$. Let $ \sigma :V(T_n)\rightarrow \{\pm 1\}$ be an assignment
of $\pm 1$ to the vertices of $T_n$ such that vertices that are
adjacent in $G$ always have opposite signs. It is easy to see that
such a $\sigma $ exists (e.g.\ let $ \sigma (v)=(-1)^{d(v_{0},v)}$
where $d(v_{0},v)$ is the number of edges in $ T_n $ that connect
$v_{0}$ to $v$). Furthermore, $\sigma $ is unique up to
multiplication by $-1$. For each $i=0,1,\cdots ,n$, let $S_{i}=(
\emptyset,c_{i})$ be the enhanced state with $s=\emptyset $,
$c_{i}$ assigns $1$ to $v_{k}$ and $x$ to $v_{j}$ for each $j\neq
i$. Let $ \varepsilon_1 =\Sigma _{i=0}^{n}\sigma (v_{k})S_{i}\in
C^{0}(T_{n})$. Let $\varepsilon_2 =(\emptyset,c)$ be the enhanced
state with $s=\emptyset $, $c$ assigns $x$ to each vertex $v_{j}$
for $j=0,\cdots ,n$. Then $ \varepsilon_1$ is a generator for
$\mathbb{Z}\{n\}$ and $\varepsilon_2$ is a generator for
$\mathbb{Z}\{n+1\}$ in $H^{0}(T_{n})$. An example is shown
in \figref{Basis for Trees}
below. 

\begin{figure}[ht!]\anchor{Basis for
Trees}
\begin{center}
\cl{\scalebox{.65}{\includegraphics{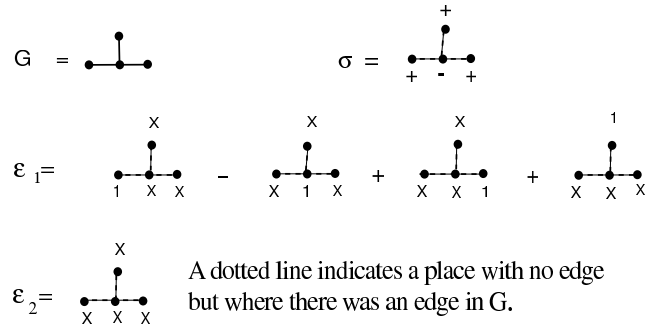}}}
\end{center}
\caption{An example of basis for trees } \label{Basis for
Trees}
\end{figure}

\begin{example}
Circuit graph with n edges \label{the circuit graph with n edges}
\end{example}

Let $G=P_{n}$ be the polygon graph with $n$ edges, also known as
the circuit graph or the cycle graph.

If $n=1$, $P_1$ is the graph with one vertex and one loop. By
Proposition \ref{graphs with loops}, $H^i(P_1)=0$ for each $i$.

 If $n=2$, $P_{2}$ is the graph with two vertices
connected by two parallel edges. By Proposition \ref{graphs with
multiple edges}, $ H^{i}(P_{2})\cong H^{i}(T_{1})\cong \left\{
\begin{array}{ll}
\mathbb{Z}\{1\}\oplus \mathbb{Z}\{2\} & \text{if }i=0 \\
0 & \text{if }i\neq 0
\end{array}
\right. $

 Next, let us assume $n>2$.

(In fact, the following method also holds for $n=2$ but the above  method is
a much easier way to get the result).

 We label the vertices of $G=P_n$ by $v_1, \cdots, v_n$
monotonically so that each $v_k$ is adjacent to $v_{k+1}$ (here
$v_{n+1}=v_1$ ). Let $e$ be the edge $v_1v_n$. Then $G-e$ is the
tree with $n$ vertices $ v_1, \cdots, v_n$ connected by a line
segment running from $v_1$ to $v_n$, and $G/e$ is the polygon
graph $P_{n-1}$ with vertices $v_1, \cdots, v_{n-1}$. The exact
sequence on $(G,e)$ gives
\begin{equation*}
\cdots \rightarrow H^{i-1}(G)\rightarrow H^{i-1}(G-e)\rightarrow
H^{i-1}(G/e)\rightarrow H^i(G)\rightarrow H^i(G-e)\rightarrow \cdots
\end{equation*}
For $i\geq 2,H^{i-1}(G-e)=H^{i}(G-e)=0$ by Example \ref{tree with
n edges}. Thus $H^{i}(G)\cong H^{i-1}(G/e)$, i.e.\
$H^{i}(P_{n})\cong H^{i-1}(P_{n-1})$ provided if $n\geq 2$ and
$i\geq 2$. Applying this equation repeatedly, we have
\begin{equation*}
H^{i}(P_{n})\cong \left\{
\begin{array}{ll}
H^1(C_{n-i+1}) & \mbox{ if } i \leq n \\
H^{i-n+1}(P_1)=0 & \mbox{ if } i \geq n.
\end{array}
\right.
\end{equation*}
Thus it suffices to determine $H^{1}(P_{n})$ and $H^{0}(P_{n})$ for
all $n$. Again, we examine part of the long exact sequence:
\begin{equation*}
0 \rightarrow H^0(G)\overset{\beta^*}{\rightarrow} H^0(G-e) \overset{
\gamma^*}{\rightarrow} H^0(G/e) \overset{\alpha^*}{\rightarrow} H^1(G)
\rightarrow 0
\end{equation*}
Here, the last group $H^1(G-e)$ is 0 because $G-e$ is a tree. This
exact sequence implies that
\begin{equation*}
H^{0}(G)\cong \ker \gamma^* ,H^{1}(G)\cong H^{0}(G/e)/\ker
\alpha^*=H^{0}(G/e)/ \mbox{Im}\gamma^*
\end{equation*}
 Thus we need to understand the map $\gamma^*
:H^{0}(G-e)\rightarrow H^{0}(G/e)$.

This map can be described as follows. An $x$ in $H^{0}(G-e)$ is
the equivalence class of a sum of terms of the form
$(\emptyset,c)$ in $ C^{0}(G-e)$. Each of these enhanced state
$S=(\emptyset ,c)$ is just a coloring of $v_{1},\cdots ,v_{n}$ by
1 or $x$. Under the map $\gamma^*$, $ S=(\emptyset ,c)$ is changed
to $(\emptyset ,\gamma (c))$ where $\gamma (c)$ is the coloring on
$V(G/e)$ defined by $\gamma (c)(v_{k})=c(v_{k})$ for each $k\neq
1, n$, and $\gamma (c)(v_{1})=c(v_{1})\ast c(v_{n})$. Basically,
for each each $(s,c)$, $\gamma ^{\ast }$ adds the edge $e$,
adjusts the colorings, then contracts $e$ to a point and
multiplies the result by $ (-1)^{i}$.

By Example (\ref{tree with n edges}), $H^{0}(G-e)\cong \mathbb{Z}
\{{n-1}\}\oplus \mathbb{Z}\{n\}$ where $\mathbb{Z}\{{n-1}\}$ is
generated by $\varepsilon _{1}$ and $\mathbb{Z}\{n\}$ is generated
by $\varepsilon _{2}$. It is easy to see that $\gamma ^{\ast
}(\varepsilon _{2})=0$ since all vertices are colored by $x$ in
$\varepsilon _2$. As for $\gamma ^{\ast
}(\varepsilon _{1})$, it will depend on the parity of $n$. We have $%
\varepsilon _{1}=S_{1}-S_{2}+\cdots +(-1)^{n-1}S_{n}$. For each
$i\neq 1,n$, $\gamma ^{\ast }(S_{i})=0$ since both $v_{1}$ and
$v_{n}$ are labeled by $x$ under $S_{i}$. For $i=1$ and $i=n $,
$\gamma ^{\ast }(S_{1})=\gamma ^{\ast }(S_{n})=\varepsilon
_{2}^{\prime } $ where $\varepsilon _{2}^{\prime }$ is the state of
$C^{0}(G/e)$ that labels every vertex by $x$. Thus $\gamma ^{\ast
}(\varepsilon _{1})=0$ if $n$ is even, and $\gamma (\varepsilon
_{2})=2\varepsilon _{2}^{\prime }$ if $n$ is odd.

It follows that $\ker$ $\gamma^* =\mbox{Span}<\varepsilon_1,
\varepsilon_2>$ if $n$ is even, and $\ker$ $\gamma^* =\mbox{Span}
<\varepsilon_2>$ if $n$ is odd. Therefore
\begin{equation*}
H^{0}(P_n)\cong \left\{
\begin{array}{ll}
\mathbb{Z}\{n\}\oplus \mathbb{Z}\{{n-1}\} & \text{ if }n \text{ is
even
and } n\geq 2 \\
\mathbb{Z}\{n\} & \text{ if }n\text{ is odd and } n>2.
\end{array}
\right.
\end{equation*}
Next, we determine $H^{1}(P_{n})$ using the same exact sequence.
We follow the discussion above. If $n$ is even, $\gamma ^{\ast
}=0$, and therefore $H^{1}(G)\cong H^{0}(G/e)\cong
\mathbb{Z}\{{n-1}\}$. If $n$ is odd, $\mbox{Im}\gamma ^{\ast
}=2\mathbb{Z}\{{n-1}\}$ in $H^{0}(G/e)$. Therefore
$H^{1}(P_{n})\cong H^{0}(G/e)/$ $\mbox{Im}\gamma ^{\ast }\cong
(\mathbb{Z}\{{n-1}\}\oplus
\mathbb{Z}\{{n-2}\})/2\mathbb{Z}\{{n-1}\}\cong
\mathbb{Z}\{{n-2}\}\oplus \mathbb{Z}_{2}\{{n-1}\}$.

As a summary, we have
\begin{equation*}
\text{For }i>0,H^{i}(P_{n})\cong \left\{
\begin{array}{ll}
\mathbb{Z}_{2}\{{n-i}\}\oplus \mathbb{Z}\{{n-i-1}\} & \text{if
}n-i\geq 2
\text{ and is even} \\
\mathbb{Z}\{{n-i}\} & \text{if }n-i\geq 2\text{ and is odd} \\
0 & \text{if }n-i\leq 1.
\end{array}
\right.
\end{equation*}
\begin{equation*}
\text{For }i=0,H^{0}(P_{n})\cong \left\{
\begin{array}{ll}
\mathbb{Z}\{n\}\oplus \mathbb{Z}\{{n-1}\} & \text{if }n\text{ is
even and }
n\geq 2 \\
\mathbb{Z}\{n\} & \text{if }n\text{ is odd and }n\geq 2 \\
0 & \text{if }n=0.
\end{array}
\right.
\end{equation*}

The following table illustrate our computational result (up to
$n=6$ and $ i=4 $) for polygon graphs.

{\def\strutt{\vrule width0pt height 11pt depth 6pt}\small\begin{center}
\begin{tabular}{|c|l|l|l|l|l|}
\hline\strutt \!\!$n \backslash i$\!\! & $H^0$ & $H^1$ & $H^2$ & $H^3$ & $H^4$ \\
\hline\strutt $P_1$ & 0 & 0 & 0 & 0 & 0 \\
\hline\strutt $P_2$ &
$\mathbb{Z}\{2\}\oplus \mathbb{Z}\{1\}$ & 0 & 0 & 0 & 0 \\
\hline \strutt$P_3$ & $\mathbb{Z}\{3\}$ & $\mathbb{Z}_2\{2\} \oplus
\mathbb{Z}\{1\}$ & 0 & 0 & 0
\\
\hline \strutt$P_4$ & $\mathbb{Z}\{4\}\oplus \mathbb{Z}\{3\}$ &
$\mathbb{Z}\{3\} $ & $\mathbb{Z }_2\{2\}\oplus \mathbb{Z}\{1\}$ &
0 &
0 \\
\hline \strutt$P_5$ & $\mathbb{Z}\{5\}$ & $\mathbb{Z}_2\{4\}\oplus
\mathbb{Z}\{3\}$ &
$\mathbb{ Z}\{3\}$ & $\mathbb{Z}_2\{2\}\oplus \mathbb{Z}\{1\}$ & 0 \\
\hline \strutt$P_6$ & $\mathbb{Z}\{6\}\oplus \mathbb{Z}\{5\}$ &
$\mathbb{Z}\{5\} $ & $\mathbb{Z }_2\{4\}\oplus \mathbb{Z}\{3\}$ &
$\mathbb{Z}\{3\}$ & $\mathbb{Z}_2\{2\}\oplus \mathbb{Z}\{1\}$\!\!
\\ \hline
\end{tabular}
\end{center}}

We note that, for all $n\geq 3$, $H^{\ast }(P_{n})$ contains
torsion. We will analyze such phenomenon for a general graph in
future work.

\section{Further developments and future problems}

After this paper was posted on the arxiv  in December 2004, there
have been some new developments. Let us mention some of these below.

First, our construction, which is based on the algebra
$A_2:=\mathbb{Z}[x]/(x^2)$, can be generalized to any commutative
graded algebras with a finite dimension at each degree. Such a
generalization has been studied in \cite{HR05}.

Some computer computations of these cohomology groups have been
done by Michael Chmutov. These computations lead to the following
observation, which has now been proved by M. Chmutov, S. Chmutov,
and Rong \cite{CCR05}. Namely, if we work on the algebra
$A_2:=\mathbb{R}[x]/(x^2)$ over a field $\mathbb{R}$ of
characteristic 0, then $H^{i, n-i}(G)\cong H^{i+1, n-i-2}(G)$ for
all $i$ and all $G$ expect for $i=0$ and $G$ is bipartite
(``knight move" isomorphism). An analogous result for Khovanov's
cohomology for alternating links/knots was proved by Lee
\cite{L02}. While this implies that the free part of the
cohomology groups is equivalent to the chromatic polynomial, it
does not say anything about torsion. Moreover, we think this is
due to fact that algebra used, $A_2$, is too simple. For more
general algebras, we do not expect such equivalence result.
Furthemore, some of the ``twisted homologies" defined in
\cite{HR05} are stronger than the chromatic polynomial.

Let us add that the isomorphism classes of chain complexes are
well defined graph invariants. They are certainly stronger than
the cohomology groups since they can distinguish some graphs with
loops.

The torsion of these homologies remains mysterious. Some studies
have been carried out on torsion in \cite{HPR05}, in joint work
with J. Przytycki. In particular, we determine precisely, when the
algebra is $A_2$, those graphs whose cohomology contains torsion.
It remains to see where tosions occur and what torsion can occur.

Finally, we ask several questions that arise naturally in our work.

\medskip
{\bf Problem 1}\qua {\sl  Understand geometric meanings of
these cohomology groups}
 \newline The chromatic polynomial has a clear geometric
interpretation. It is not clear what our cohomology groups measure.
We note that these cohomology groups are not invariants of
  matroid type of the graphs. For example, the graph made of
two triangles glued at one vertex and the graph which is the
disjoint union of two triangles have the same matroid type but
different chromatic polynomials, and therefore different
cohomology groups.

\medskip
{\bf Problem 2}\qua {\sl Investigate functorial properties of
the cohomology groups}
 \newline  The classical homology theory is a functor:
continuous maps between spaces induce homomorphism between
homology groups.  Khovanov's link cohomology $Kh$ also satisfies a
functorial property: a cobordism $C$ between two links $L_1$ and
$L_2$ induces a homomorphism between the Khovanov homologies of
$L_1$ and $L_2$.
 For our graph cohomology, we can associate homomorphisms between
 cohomology groups to each inclusion map of graphs. Essentially, it
 is the iteration of the map $\beta^*$ in the long exact
 sequence. We are not sure if this could lead to obstructions of one graph
being a subgraph of another. In particular, can this be used to
detect $K_5$ and $K_{3,3}$ in a given graph, and therefore
determine the planarity of the graph? In any case, further study
of this functor is worthwhile.

\medskip
{\bf Problem 3}\qua {\sl  Carry out further study of torsions}
\newline  First, what orders of torsion can we have? This will
reflect the algebraic properties of the underlying algebra.  As
algebra varies, we can get torsions of arbitrary orders
\cite{HR05} \cite{HPR05}. For the algebra $A_2$, the only torsions
based on examples we know are of order two. We would like to know
whether other orders can occur. We would also like to know whether
the distributions of torsion have any pattern, e.g.\ some kind of
modified  knight moves.

\medskip
{\bf Problem 4}\qua {\sl Study connections with Khovanov's
link cohomology and other homologies}  \newline It is natural to
ask for connections with link homologies given the connections
between graph polynomials and link polynomials. A direct relation
between $H^i(G)$ and the framed link cohomology (which corresponds
to the Kauffman bracket) was given in \cite{HPR05}, for all $i$ up
to the length of the shortest cycle in $G$. Also, J. Przytycki has
pointed out that our work in [HPR05] also shows connections of graph
cohomology with other homology theories such as Hochschild homology
of algebras.

\medskip
{\bf Problem 5}\qua {\sl Establish cohomology theories for
various other graph polynomials, and for signed graphs}
\newline This is
a natural question. Some work has already been done in this
direction following our current paper. See, e.g.\
\cite{S05a}\cite{S05b}\cite{JR05}.

\Addresses\recd
\end{document}